# Limits of tangent spaces along a subvariety at an isolated singularity

Achim Hennings [♦]

**Abstract**

We show that the dimension of the set of limits of tangent spaces along a subvariety at an irreducible isolated singularity is one less than the dimension of the subvariety itself.

In the note [4] a proof was given, that the set of limits of tangent spaces at an irreducible isolated singularity has dimension one less than the dimension of the singularity, thus answering a question posed by A. Simis, K. Smith and B. Ulrich [6] and also by B. Youssin [8]. In the present note we prove a similar statement for limits along a subvariety. We employ a different method, which is close to the arguments of W. Fulton and J. Hansen [3] and of F.L. Zak [9] in the projective case.

Let $(X, 0) \subseteq (\mathbb{C}^N, 0)$ be a reduced isolated singularity of dimension $d \geq 1$. To exclude the trivial case of smooth components only, we assume at once that $(X, 0)$ is irreducible and non-smooth. Let $X^*$ be the smooth locus of a sufficiently small representative $X$.

We recall the geometric construction of the Nash transformation ([5], [7]). Let $\gamma: X^* \to G$, $\gamma(x) = T_{X,x}$, be the Gauss map with values in the Grassmann manifold of $d$−dimensional subspaces of $\mathbb{C}^N$, assigning tangent spaces to smooth points. The Nash transform is the closure $\hat{X} \subseteq X \times G$ of the graph of $\gamma$. The restriction of the first projection gives a proper modification

$$\nu: \hat{X} \to X,$$

also called Nash transformation. The exceptional set is the special fibre $\nu^{-1}(0) \subseteq G$, consisting of all limits of tangent spaces at points of $X^*$ which converge to $0$.

Now we take a subvariety (irreducible germ of analytic set) $(Y, 0) \subseteq (X, 0)$ of dimension $\geq 1$ with representative $Y$ and put $Y^* = X^* \cap Y$. It is reasonable to call the closure $\hat{Y} = \overline{Y^*}$ in $\hat{X}$ the strict transform of $Y$. The special fibre of $\hat{Y}$, denoted $(\hat{Y})_0$, is the set of limits of tangent spaces of $X$ along $Y$, that is, we are considering them only at points of $Y^*$ instead of $X^*$. Our result is:

**Theorem:** The dimension of $(\hat{Y})_0$ attains the largest possible value $\dim(Y, 0) - 1$.

We note that in case $\dim Y = 1$ the assertion is true. For $\dim Y = 2$ it is contained in the following

**Lemma:** If $\dim(Y, 0) \geq 2$, then $\dim(\hat{Y})_0 \geq 1$.

Proof: We assume on the contrary that $(\hat{Y})_0$ is finite. Then a general linear projection $\pi: \mathbb{C}^N \to \mathbb{C}^d$ is injective when restricted to the finitely many limit spaces $T \in (\hat{Y})_0$. Upon choosing a smaller representative, the restriction $\pi_X = \pi|X$ is then unramified along $Y^*$. We consider the map

[♦] Universität Siegen, Fakultät IV, Hölderlinstrasse 3, D-57068 Siegen

$$f = (\pi|X) \times (\pi|Y) : X \times Y \to \mathbb{C}^d \times \mathbb{C}^d.$$

Let $\Delta_Y \subseteq X \times Y$ and $\Delta \subseteq \mathbb{C}^d \times \mathbb{C}^d$ be the diagonal subspaces. The local injectivity of $\pi_X$ along $Y^*$ implies that $\Delta_Y \cap X^* \times Y^*$ is open in $f^{-1}(\Delta) \cap X^* \times Y^*$. (Otherwise there would be sequences $x_n \in X^*$, $y_n \in Y^*$, with $x_n \neq y_n$, $\pi_X(x_n) = \pi_X(y_n)$ and a common limit $y \in Y^*$, a contradiction.) We have $f^{-1}(\Delta) \cap X^* \times Y^* = f^{-1}(\Delta) - \{(0,0)\}$, and $\Delta_Y \cap X^* \times Y^* = \Delta_Y - \{(0,0)\}$ is closed therein. Since $f^{-1}(\Delta)$ is cut out of $X \times Y$ by $d \leq \dim X \times Y - 2$ equations, the germ of $f^{-1}(\Delta) - \{(0,0)\}$ at $(0,0)$ is connected by a theorem of Grothendieck ([1], [2]). It follows that, possibly for a smaller representative $X$, $f^{-1}(\Delta) - \{(0,0)\} = \Delta_Y - \{(0,0)\}$. Hence $\pi_X$ is injective on the whole preimage $\pi_X^{-1}(\pi_X(Y))$ and thus has degree 1. But this could only be true if $X$ was smooth.

Proof of the theorem: For the proof of the theorem we use a representation of the Nash transform as a subspace of projective space over $X$. By the Plücker embedding the Grassmann manifold $G = G(d, N)$ can be viewed as a subspace of $\mathbb{P}^m$, where $m = \binom{N}{d} - 1$. This gives an embedding of the Nash transform $\hat{X}$ into $X \times \mathbb{P}^m$. For a hyperplane $H \subseteq \mathbb{P}^m$ we can take the intersection of $\hat{X}$ with the divisor $\bar{H} = X \times H$ of $X \times \mathbb{P}^m$. If $\dim (\hat{X})_0 > 0$ and $H$ intersects the irreducible components of $(\hat{X})_0$ properly, then $\hat{X} \cap \bar{H}$ has an irreducible component near $(\hat{X})_0$ of codimension 1, not contained in $(\hat{X})_0$. The image of this component under projection to $(X, 0)$ is also a subvariety of codimension 1.

These arguments can be applied in the same way to the strict transform $\hat{Y}$ of a subvariety $(Y, 0)$ of $(X, 0)$. Hence we obtain a subvariety $(Z, 0) \subseteq (Y, 0)$ of codimension 1 such that $\dim(\hat{Z})_0 < \dim(\hat{Y})_0$, provided that $\dim(\hat{Y})_0 \geq 1$.

Now we assume that there is a subvariety $(Y, 0)$ of $(X, 0)$ of dimension $\geq 2$ having $\dim(\hat{Y})_0 < \dim Y - 1$. We choose $Y$ of minimal dimension. By the lemma, $\dim Y \geq 3$ and $\dim(\hat{Y})_0 \geq 1$. By the above consideration there is a subvariety $(Z, 0) \subseteq (Y, 0)$ of codimension 1 such that $\dim(\hat{Z})_0 < \dim(\hat{Y})_0$. This contradicts the choice of $Y$.